\documentclass[11pt,a4paper]{article}
\usepackage[top=1in, bottom=1.25in, left=1in, right=1in]{geometry}
\usepackage[utf8]{inputenc}
\usepackage[english]{babel}
\usepackage{graphicx,epstopdf}
\usepackage{caption, subcaption,float}
\usepackage{color}
\usepackage{multirow}
\usepackage{tikz}
\usepackage{amsthm}
\usepackage{tabls}
\usepackage{cite}
\usepackage[round,comma,sort]{natbib}

\usepackage[ruled, linesnumberedhidden,vlined]{algorithm2e} 
\usetikzlibrary{matrix}
\usepackage{amsmath,amssymb,dsfont,mathtools}
\usepackage{hyperref}
\hypersetup{
    colorlinks=true,       
    linkcolor=blue,          
    citecolor=olive,       
    filecolor=magenta,      
    urlcolor=blue          
}
\usepackage{aliascnt}

\usepackage{soul}

\begin{document}

\title{Intersection of Parabolic Subgroups in \\ Euclidean Braid Groups: a short proof}
\date{\today }
\author{Mar\'{i}a Cumplido, Federica Gavazzi and Luis Paris}


\maketitle
\theoremstyle{plain}
\newtheorem{theorem}{Theorem}

\newaliascnt{lemma}{theorem}
\newtheorem{lemma}[lemma]{Lemma}
\aliascntresetthe{lemma}
\providecommand*{\lemmaautorefname}{Lemma}

\newaliascnt{proposition}{theorem}
\newtheorem{proposition}[proposition]{Proposition}
\aliascntresetthe{proposition}
\providecommand*{\propositionautorefname}{Proposition}

\newaliascnt{corollary}{theorem}
\newtheorem{corollary}[corollary]{Corollary}
\aliascntresetthe{corollary}
\providecommand*{\corollaryautorefname}{Corollary}

\newaliascnt{conjecture}{theorem}
\newtheorem{conjecture}[conjecture]{Conjecture}
\aliascntresetthe{conjecture}
\providecommand*{\conjectureautorefname}{Conjecture}

\newaliascnt{question}{theorem}
\newtheorem{question}[question]{Question}
\aliascntresetthe{question}
\providecommand*{\conjectureautorefname}{Question}

\theoremstyle{remark}

\newaliascnt{claim}{theorem}
\newaliascnt{remark}{theorem}

\newtheorem{claim}[claim]{Claim}
\newtheorem{remark}[remark]{Remark}
\newaliascnt{notation}{theorem}
\newtheorem{notation}[notation]{Notation}
\aliascntresetthe{notation}
\providecommand*{\notationautorefname}{Notation}

\aliascntresetthe{claim}
\providecommand*{\claimautorefname}{Claim}

\aliascntresetthe{remark}
\providecommand*{\remarkautorefname}{Remark}

\newtheorem*{claim*}{Claim}
\theoremstyle{definition}

\newaliascnt{definition}{theorem}
\newtheorem{definition}[definition]{Definition}
\aliascntresetthe{definition}
\providecommand*{\definitionautorefname}{Definition}

\newaliascnt{example}{theorem}
\newtheorem{example}[example]{Example}
\aliascntresetthe{example}
\providecommand*{\exampleautorefname}{Example}

\def\autorefspace{\hspace*{-0.5pt}}
\def\sectionautorefname{Section\autorefspace}
\def\subsectionautorefname{Section\autorefspace}
\def\subsubsectionautorefname{Section\autorefspace}
\def\figureautorefname{Figure\autorefspace}
\def\subfigureautorefname{Figure\autorefspace}
\def\tableautorefname{Table\autorefspace}
\def\equationautorefname{Equation\autorefspace}
\def\Itemautorefname{item\autorefspace}
\def\Hfootnoteautorefname{footnote\autorefspace}
\def\AMSautorefname{Equation\autorefspace}

\newcommand{\co}{\simeq_c}
\newcommand{\w}{\widetilde}
\newcommand{\po}{\preccurlyeq}
\newcommand{\dist}{\mathrm{d}}

\def\Z{\mathbb Z} 
\def\Ker{{\rm Ker}} \def\R{\mathbb R} \def\GL{{\rm GL}}
\def\HH{\mathcal H} \def\C{\mathbb C} \def\P{\mathbb P}
\def\SSS{\mathfrak S} \def\BB{\mathcal B} \def\PP{\mathcal P} 
\def\supp{{\rm supp}} \def\Id{{\rm Id}} \def\Im{{\rm Im}}
\def\MM{\mathcal M} \def\S{\mathbb S}
\newcommand{\bigveer}{\bigvee^\Lsh}
\newcommand{\wedger}{\wedge^\Lsh}
\newcommand{\veer}{\vee^\Lsh}
\def\diam{{\rm diam}}

\newcommand{\myref}[2]{\hyperref[#1]{#2~\ref*{#1}}}

\begin{abstract}
We give a short proof for the fact, already proven by Thomas Haettel, that the arbitrary intersection of parabolic subgroups in Euclidean Braid groups $A[\tilde{A}_n]$ is again a parabolic subgroup. To that end, we use that the spherical-type Artin group $A[B_{n+1}]$ is isomorphic to  $A[\tilde{A}_n]\rtimes \Z$ . 

\medskip

{\footnotesize
\noindent \emph{2020 Mathematics Subject Classification.} 20F36.

\noindent \emph{Key words.} Artin groups; Euclidean braid groups; parabolic subgroups.}

\end{abstract}

An Artin (or Artin-Tits) group $A_S$ is any group defined by a finite set of generators $S$ and relations of the type $sts\cdots=tst\cdots$, for $s,t\in S$, where each word in the equality has the same length, denoted $m_{s,t}$. The number of known global results for these groups is very limited, and for some decades now, classic problems such as the word problem, the conjugacy problem, or the $K(\pi, 1)$ conjecture have been the subject of study by group theorists. Specifically, it has become necessary to better study the properties of certain specific subgroups: the parabolic subgroups. A standard parabolic subgroup $A_X$ of $A_S$ is a subgroup generated by a subset $X\subset S$, and thanks to \citep{Vanderlek}, we also know that it coincides with the Artin group on $X$ with the same relations that these generators have in $A_S$. A parabolic subgroup is any conjugate of a standard parabolic subgroup. These subgroups play a principal role in the construction of simplicial complexes associated with Artin groups \citep{CharneyDavis,CGGW,CMV,ConciniSalvetti,LP12kp1}, either as stabilizers of simplices or as the building blocks of the complexes. In the case of braid groups, the parabolic subgroups coincide with the isotopy classes of non-degenerate multicurves in the $n$-punctured disc (containing the vertices of the same curve complex). However, the basic question of whether the arbitrary intersection of parabolic subgroups is a parabolic subgroup remains open in most cases.

To establish notations, we will use Coxeter graphs. A Coxeter graph $\Gamma_S$ encodes the information of an Artin group $A_S$ as follows: each generator corresponds to a vertex, and two vertices $s,t$ are connected by an edge if the vertices do not commute. This edge is labeled by $m_{s,t}$ if $m_{s,t}>3$ and by $\infty$ if there is no relation between $s$ and $t$. In this way, we can also refer to $A_S$ as $A[\Gamma_S]$. The Euclidean braid group (also known as affine braid group) with $n+1$ generators is the group $A[\tilde{A}_n]$, where $\tilde{A}_n$ is the graph in  \autoref{graphtildean}:
\begin{figure}[h]
    \centering
    \includegraphics[width=4.5cm]{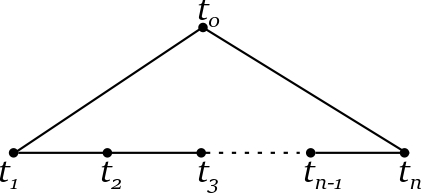}
    \caption{The Coxeter graph $\Tilde{A}_n$.}
    \label{graphtildean}
\end{figure}

 In this article, we will give an alternative and concise proof of the following theorem:

\begin{theorem}[{\citealp[Corollary~N]{Haettel}}]\label{maintheorem}
The arbitrary intersection of parabolic subgroups in~$A[\tilde{A}_n]$ is a parabolic subgroup.    
\end{theorem}

\begin{figure}[h]
    \centering
    \includegraphics[width=5.5cm]{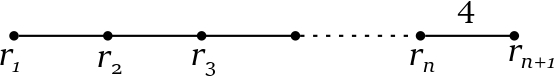}
    \caption{The Coxeter graph $B_{n+1}$.}
    \label{graphBn}
\end{figure}

We know, according to \citep{CGGW}, that the intersection of parabolic subgroups of a spherical-type Artin group remains a parabolic subgroup. In particular, this is true in the group $A[B_{n+1}]$, where $B_{n+1}$ is the graph illustrated in \autoref{graphBn}.

Using the notation of the figure, let $\rho=r_1\dots r_n r_{n+1}$. For $1\leq i \leq n-1$, it holds that $\rho r_i \rho^{-1}=r_{i+1}$. If we additionally define $r_0:= \rho r_n \rho^{-1}$, we can state that the previous equality is true modulo $n + 1$ (observing that $\rho^2 r_{n}\rho^{-2} = r_1$ is sufficient). Furthermore, there exists an exterior automorphism $f$ of $A[\tilde{A}_n]$ that sends $t_i$ to $t_{i+1}$ modulo $n+1$, and an action of $\Z\cong \langle u \rangle$ on  $A[\tilde{A}_n]$ defined by setting $u^{-1} g u=f (g)$. With this action, we can define the semidirect product $A[\tilde{A}_n] \rtimes \langle u \rangle$.

\begin{theorem}{\citep{KP02}}
The map $\varphi: A[\tilde{A}_n] \rtimes \langle u \rangle \longrightarrow A[{B}_{n+1}] $ that sends $t_i$ to~$r_i$ and $u$ to $\rho$ is an isomorphism.
    
\end{theorem}

\begin{remark}
The restriction of $\varphi$ to $A[\tilde{A}_n]$ gives an embedding of this Artin group in $A[B_{n+1}]$.
Consider the group homomorphism $\xi: A[B_{n+1}]\longrightarrow \Z$ defined by  $r_i\longmapsto 0$ for $1\leq i \leq n$ and $r_{n+1}\longmapsto 1$. Observe that $\rho$ is mapped to $1$, and that the kernel of $\xi$ is $\varphi(A[\tilde{A}_n])$.
\end{remark}
\begin{proof}[Proof of \autoref{maintheorem}]
    Firstly, notice that if $P$ is a proper parabolic subgroup of $A[\tilde{A}_n]$, then~$\varphi(P)$ is a parabolic subgroup of $A[B_{n+1}]$. The only case in which this is not clear is when $P=gA[\tilde{A}_n]_Xg^{-1}$ and $X$ contains $t_0$, which is sent to $\rho r_n \rho^{-1}$ by $\varphi$. However, since $P$ is proper, $X$ does not contain all the generators of $A[\tilde{A}_n]$, thus, using $\rho$, we can always conjugate $\varphi(X)$ to a subset of $\{r_1,\dots, r_n\}$.\\
    Now suppose that $P_1$ and $P_2$ are two parabolic subgroups of $A[\tilde{A}_n]$. Since in $A[B_{n+1}]$ the intersection of parabolic subgroups is a parabolic subgroup, we have, in particular, that $\varphi(P_1\cap P_2)=\varphi(P_1)\cap \varphi(P_2)$ is a parabolic subgroup of $A[B_{n+1}]$.  
    
    To complete the proof, it remains to show that if $Q$ is a parabolic subgroup of $A[B_{n+1}]$ such that $P:=\varphi^{-1}(Q)\subset A[\tilde{A}_n]$, then $P$ is a parabolic subgroup of $A[\tilde{A}_n]$. We can write $Q=h A[B_{n+1}]_Y h^{-1}$, with $Y\subset\{r_1,\cdots,r_{n+1}\}$ and $h\in A[B_{n+1}]$. First, we show that $r_{n+1}\not\in Y$. If we suppose otherwise, then $hr_{n+1}h^{-1}\in Q$ and $\varphi^{-1}(hr_{n+1}h^{-1})\in P\subset A[\tilde{A}_n]$. In particular $hr_{n+1}h^{-1}$ must belong to the kernel of $\xi:A[B_{n+1}]\longrightarrow \Z$, that coincides with $\varphi(A[\Tilde{A}_n])$. But computing $\xi(hr_{n+1}h^{-1})=1$, we get a contradiction. Therefore, $Y\subset\{r_1,\dots,r_n\}= \varphi(\{t_1,\dots,t_n\})$. Furthermore, since  $A[\tilde{A}_n] \rtimes \langle u \rangle \cong A[{B}_{n+1}]$, we can write $h=h_1\rho^m$, with $\varphi^{-1}(h_1)\in A[\tilde{A}_n]$ and $m\in \Z$. Thus, $Q=h_1\rho^m A[B_{n+1}]_Y \rho^{-m} h_1^{-1}$, and exploiting the fact that $u^m A[\Tilde{A}_n]_{\varphi^{-1}(Y)} u^{-m}=f^m(A[\Tilde{A}_n]_{\varphi^{-1}(Y)})$, we get
    $$ P=\varphi^{-1}(Q) = \varphi^{-1}(h_1) f^m (A[\tilde{A}_n]_{\varphi^{-1}(Y)}) (\varphi^{-1}(h_1))^{-1}= \varphi^{-1}(h_1) A[\tilde{A}_n]_{f^m (\varphi^{-1}(Y))} (\varphi^{-1}(h_1))^{-1},$$
    which has the form of a parabolic subgroup in $A[\Tilde{A}_n]$.
    Finally ---see the complete argument in \cite[Corollary~16]{CMV}---, in any Artin group, any descending chain of inclusions of parabolic subgroups has to stabilize, and then the finite intersection of parabolic subgroups, being a parabolic subgroup, implies that the same is true for an arbitrary intersection.  
\end{proof}

\medskip
\newpage

\noindent{\textbf{\Large{Acknowledgments}}}

The first author was supported by a Ram\'on y Cajal 2021 grant and the research grant PID2022-138719NA-I00 (Proyectos de Generaci\'on de Conocimiento 2022), both financed by the Spanish Ministry of Science and Innovation.

\medskip
\bibliography{Bibliography}

\begin{thebibliography}{}

\bibitem[\protect\citename{Charney \& Davis, }1995]{CharneyDavis}
Charney, Ruth, \& Davis, Michael~W. 1995.
\newblock The $K(\pi, 1)$-Problem for Hyperplane Complements Associated to
  Infinite Reflection Groups.
\newblock {\em J. Amer. Math. Soc.}, {\bf 8}(3), 597--627.

\bibitem[\protect\citename{Concini \& Salvetti, }2000]{ConciniSalvetti}
Concini, Corrado, \& Salvetti, Mario. 2000.
\newblock Cohomology of Coxeter groups and Artin groups.
\newblock {\em Mathematical Research Letters}, {\bf 7}(03), 213--232.

\bibitem[\protect\citename{Cumplido {\em et~al.}, }2019]{CGGW}
Cumplido, Mar\'ia, Gebhardt, Volker, González-Meneses, Juan, \& Wiest, Bert.
  2019.
\newblock On parabolic subgroups of {Artin–Tits} groups of spherical type.
\newblock {\em Adv. Math.}, {\bf 352}, 572 -- 610.

\bibitem[\protect\citename{Cumplido {\em et~al.}, }2023]{CMV}
Cumplido, Mar\'ia, Martin, Alexandre, \& Vaskou, Nicolas. 2023.
\newblock {Parabolic subgroups of large-type Artin groups}.
\newblock {\em Mathematical Proceedings of the Cambridge Philosophical
  Society}, {\bf 174}(2), 393–414.

\bibitem[\protect\citename{Haettel, }2021]{Haettel}
Haettel, Thomas. 2021.
\newblock {\em Lattices, injective metrics and the $K(\pi,1)$ conjecture}.
\newblock Algebraic \& Geometric Topology (to appear). arXiv:2109.07891.

\bibitem[\protect\citename{Kent~IV \& Peifer, }2002]{KP02}
Kent~IV, Richard, \& Peifer, David. 2002.
\newblock A Geometric and Algebraic Description of Annular Braid Groups.
\newblock {\em IJAC}, {\bf 12}(02), 85--97.

\bibitem[\protect\citename{Paris, }2014]{LP12kp1}
Paris, Luis. 2014.
\newblock $K(\pi ,1)$ conjecture for {Artin} groups.
\newblock {\em Annales de la Facult\'e des sciences de Toulouse :
  Math\'ematiques}, {\bf Ser. 6, 23}(2), 361--415.

\bibitem[\protect\citename{Van~der Lek, }1983]{Vanderlek}
Van~der Lek, Harm. 1983.
\newblock {\em {The Homotopy Type of Complex Hyperplane Complements}}.
\newblock Ph.D. thesis, Nijmegen.

\end{thebibliography}

\bigskip\bigskip{\footnotesize%

\noindent
\textit{\textbf{María Cumplido} \\ 
Departamento de \'Algebra,
Facultad de Matem\'aticas,
Universidad de Sevilla. \\
Calle Tarfia s/n
41012, Seville, Spain.} \par
 \textit{E-mail address:} \texttt{\href{mailto:cumplido@us.es}{cumplido@us.es}}

\noindent
\textit{\textbf{Federica Gavazzi} \\  IMB, UMR5584, CNRS, Université de Bourgogne, 21000 Dijon,France} \par
 \textit{E-mail address:} \texttt{\href{mailto:Federica.Gavazzi@u-bourgogne.fr}{Federica.Gavazzi@u-bourgogne.fr}}

\noindent
\textit{\textbf{Luis Paris} \\  IMB, UMR5584, CNRS, Université de Bourgogne, 21000 Dijon,France} \par
 \textit{E-mail address:} \texttt{\href{mailto:luis.paris@u-bourgogne.fr}{luis.paris@u-bourgogne.fr}}
 }

\end{document}